\documentclass[12pt,letterpaper]{amsart}
\usepackage[utf8]{inputenc}
\usepackage{mathtools}
\usepackage{booktabs}
\usepackage{amssymb, amsmath, enumerate, amsthm, cite, colonequals, stmaryrd, lipsum, graphicx, fancyvrb,fancyhdr, indentfirst, times, setspace}
\usepackage[margin=1in]{geometry}
\usepackage{mathtools}
\makeatletter
\@namedef{subjclassname@2020}{\textup{2020} Mathematics Subject Classification}
\makeatother

\usepackage[colorlinks]{hyperref}
\hypersetup{
 colorlinks=true,
 linkcolor=blue,
 filecolor=magenta, 
 urlcolor=cyan,
}

\usepackage[capitalise]{cleveref}
\crefformat{equation}{(#2#1#3)}
\crefrangeformat{equation}{(#3#1#4--#5#2#6)}
\crefformat{enumi}{(#2#1#3)}
\crefrangeformat{enumi}{(#3#1#4--#5#2#6)}

\usepackage[pagewise]{lineno}
\let\oldequation\equation
\let\oldendequation\endequation
\renewenvironment{equation}{\linenomathNonumbers\oldequation}{\oldendequation\endlinenomath}
\expandafter\let\expandafter\oldequationstar\csname equation*\endcsname
\expandafter\let\expandafter\oldendequationstar\csname endequation*\endcsname
\renewenvironment{equation*}{\linenomathNonumbers\oldequationstar}{\oldendequationstar\endlinenomath}
\let\oldalign\align
\let\oldendalign\endalign
\renewenvironment{align}{\linenomathNonumbers\oldalign}{\oldendalign\endlinenomath}
\expandafter\let\expandafter\oldalignstar\csname align*\endcsname
\expandafter\let\expandafter\oldendalignstar\csname endalign*\endcsname
\renewenvironment{align*}{\linenomathNonumbers\oldalignstar}{\oldendalignstar\endlinenomath}

\usepackage{tikz-cd}




\newcounter{intro}
\newcounter{result}

\newtheorem{introthm}[intro]{Theorem}

\theoremstyle{definition}
\newtheorem{theorem}[result]{Theorem}
\newtheorem{lemma}[result]{Lemma}

\newtheorem{proposition}[result]{Proposition}
\newtheorem{corollary}[result]{Corollary}

\newtheorem{remark}[result]{Remark}

\numberwithin{equation}{subsection}

\title[Stability conditions on Calabi--Yau threefolds of type K]{Stability conditions on Calabi--Yau threefolds of type K}

\author[U.~Dutta]{Uttaran Dutta}
\address{Department of Mathematics, University of South Carolina}
\email{udutta@email.sc.edu}

\date{\today}

\definecolor{bleudefrance}{rgb}{0.19, 0.55, 0.91}

\def\P{\mathbb{P}}
\def\C{\mathbb{C}}

\def\Q{\mathbb{Q}}

\def\L{\mathbb{L}}

\def\Z{\mathbb{Z}}

\def\DD{\mathcal{D}}

\def\TT{\mathcal{T}}
\def\PP{\mathcal{P}}

\def\R{\mathbb{R}}
\def\OO{\mathcal{O}}
\def\FF{\mathcal{F}}
\def\NN{\mathcal{N}}
\def\AA{\mathcal{A}}

\def\OO{\mathcal{O}}

\DeclareMathOperator{\ch}{ch}
\DeclareMathOperator{\Db}{\mathrm{D}^b}

\DeclareMathOperator{\Aut}{Aut}
\DeclareMathOperator{\rk}{rk}

\DeclareMathOperator{\Hom}{Hom}

\DeclareMathOperator{\Pic}{Pic}

\DeclareMathOperator{\NS}{NS}

\DeclareMathOperator{\Supp}{Supp}

\DeclareMathOperator{\Coh}{Coh}
\DeclareMathOperator{\sHom}{\mathcal{H}\kern -.5pt\mathit{om}}
\DeclareMathOperator{\sTor}{\mathcal{T}\kern -1.5pt\mathit{or}}
\DeclareMathOperator{\GL}{GL}

\DeclareMathOperator{\Stab}{Stab}

\DeclareMathOperator{\td}{td}

\DeclareMathOperator{\Knum}{K_{\mathrm{num}}}

\begin{document}

\maketitle
\setcounter{tocdepth}{1}
\numberwithin{result}{section}
\begin{abstract}
    A Calabi--Yau threefold of type K is a free quotient $X=(S\times E)/G$ of the product of a K3 surface and an elliptic curve. We construct a family of geometric Bridgeland stability conditions on $X$ by using Liu's product construction on $S\times E$ and then inducing along the Galois covering. The family has dimension $2h^{1,1}(X)$, and is disjoint from the $\widetilde{\GL}_2^+(\R)$-orbit of Li's stability conditions. These threefolds carry a K3 fibration, and the pushforward of an object from a smooth K3 fiber is (semi)stable if and only if it is (semi)stable on the K3 surface.
\end{abstract}
\section{Introduction}
Let $X$ be a Calabi--Yau threefold over $\C$ with infinite fundamental group. By the Bogomolov decomposition theorem, the minimal splitting covering of $X$ is either an abelian threefold (\emph{type A}) or a product $S\times E$ (\emph{type K}), where $S$ is a K3 surface and $E$ an elliptic curve. Oguiso--Sakurai \cite{Oguiso-sakurai-type-A} classified type A threefolds and Hashimoto--Kanazawa \cite{TypeK-classification} classified type K threefolds. \cite{TypeK-classification} classification divided type K threefolds into eight deformation equivalence classes, distinguished by the Galois group (See Theorem \ref{thm:typeK-classification}).

Bridgeland stability conditions exist on every smooth projective varieties following \cite{Li}, who restricted stability conditions on $\P^n$ along a closed embedding. Recently, \cite{LLLMPSZ:26} extended this to projective families and \cite{Cheng:26} proved that these stability conditions satisfy the support property with respect to the full numerical Grothendieck group. \cite{Cheng-Feyzbaksh} proved that on a threefold, for large volume, the heart of Li's stability conditions is the double-tilt heart of Bayer--Macr\`i--Toda \cite{BMT-I}. Previously, a Bogomolov--Gieseker type inequality was proved for abelian threefolds \cite{BMS16}, but the inequality fails for a Calabi--Yau threefold with infinite fundamental group \cite{FKLR-Brill-Noether-threefolds}. Very recently, \cite{LLPZ-Frobenius} have established a Bogomolov--Gieseker-type inequality for tilt-semistable objects on smooth projective varieties of dimension $\geq 3.$

Following \cite{Liu-product}, we first construct stability conditions on $S\times E$ and then induce them to $X$ \cite{MMS,Perr-Shah}. As mentioned in \cite[Example 3.7]{FKLR-Brill-Noether-threefolds}, the existence of stability conditions on type K threefolds follows in this way from \cite{Liu-product} and \cite{Perr-Shah}.

We make this construction explicit, obtaining stability conditions $\sigma^G_{\omega,B,s,t}$ on $X$ for $\omega,B\in\NS(S)^G_{\Q}$ and $s,t\in\Q_{>0}.$ In \cite{FKLR-Brill-Noether-threefolds}, it is noted that the stability conditions need not be geometric in general (cf.\ \cite[Remark 1.3]{Perr-Shah}). We show that the stability conditions $\sigma^G_{\omega,B,s,t}$ constructed here are geometric: every skyscraper sheaf is $\sigma^G_{\omega,B,s,t}$-stable of the same phase. 

\subsection*{Results}
Let $X=(S\times E)/G$ be a type K threefold and $\pi:S\times E\to X$ be the covering. Let $\omega,B\in\NS(S)^G$ with $\omega$ ample and $\omega^2>2.$ Following \cite{Bri-K3}, set $\sigma_{\omega,B}$ for the stability conditions on $S.$ For $s,t\in\Q_{>0},$ Liu's construction gives us a stability condition $\sigma_{\omega,B,s,t}$ on $S\times E.$

\begin{introthm}
    (Theorem \ref{mainthm}, Corollary \ref{cor:main2}) $\sigma_{\omega,B,s,t}$ is $G$-invariant and induces a stability condition on $X.$ An object $F\in\Db(X)$ is $\sigma^G$-semistable if and only if $\pi^*F$ is $\sigma_{\omega,B,s,t}$-semistable. All skyscraper sheaves are stable with phase one. The family has dimension $2\rk\NS(S)^G+2=2h^{1,1}(X)$.
\end{introthm}

Let $Y$ be a smooth projective variety with a rational numerical stability condition $\sigma.$ Let $E$ be an elliptic curve and $\sigma^{L,s,t}$ be Liu's stability condition on $Y\times E,$ where $s,t,\in\Q_{>0}.$ For $e\in E$, let $\DD_e$ be the subcategory of $\Db(Y\times E)$ of objects supported on $Y\times\{e\}.$

\begin{introthm}
    (Theorem \ref{mainthm: stability-fiber}, \ref{cor:main1},) An object $F\in\Db(Y)$ is a $\sigma$-semistable (resp. stable) if and only if $i_{e*}F$ is $\sigma^{L,s,t}$-semistable (resp. stable). For a type K threefold this holds for objects supported on a smooth K3 fiber of the K3-fibration $X\to\P^1$.
\end{introthm}

\begin{introthm}
    (Theorem \ref{thm:Li-vs-Liu}) For every closed embedding $X\hookrightarrow\P^n$ and all parameters, the stability conditions $\sigma^G_{\omega,B,s,t}$ do not lie in the $\widetilde{\GL}_2^+(\R)$-orbit of Li's stability conditions.
\end{introthm}

\subsection*{Conventions}
We work over $\C$. For a smooth projective variety $Y$, $\Db(Y)$ is the bounded derived category of coherent sheaves and $\Knum(Y)$ is its numerical Grothendieck group. A stability condition $\sigma=(\AA,Z)$ is with respect to a pair $(\Lambda,v)$ and satisfies the support property. It is \emph{numerical} if $v$ factors through $\Knum(Y)$,

\subsection*{Acknowledgement}
 The author thanks Nicol\'as Vilches, conversations with whom on a separate project led the author to \cite[\S 3]{LMPSZ}, upon which Section \ref{sec:fiber} is built. The author was partially supported by a Simons Dissertation Fellowship in Mathematics.

\subsection*{AI disclosure}
The author used Anthropic's Claude Opus 5.5 to proofread the manuscript, which found many typographical errors and structural inconsistencies. It has not, however, produced any of the statements or proofs. All of the contents of this article were conceived of and written by the author.

\section{Preliminaries}
\subsection{Calabi--Yau threefolds of type K}
The Bogomolov decomposition theorem implies that a Calabi–Yau threefold $X$ with infinite
fundamental group admits an \'etale Galois covering either by an abelian threefold or by the product of a K3 surface and an elliptic curve. The threefold $X$ is called of type K in the latter case \cite{TypeK-classification}. Among many candidates of such coverings, we can always find a unique
smallest one, up to isomorphism as a covering. We call the smallest covering the minimal
splitting covering of $X$.

\begin{theorem}\label{thm:typeK-classification}
    \cite[Theorem 3.1]{TypeK-classification} There exist exactly eight Calabi–Yau threefolds of type K, up to deformation equivalence. The equivalence class is uniquely determined by the Galois group $G$ of the minimal splitting covering. Moreover, the Galois group is isomorphic to one of the following combinations of cyclic and dihedral groups
    \[
    C_2, C_2\times C_2, C_2\times C_2\times C_2, D_6,D_8, D_{10}, D_{12},\text{ or }C_2\times D_8.
    \]
\end{theorem}

\begin{proposition}\label{prop:HK1}
    \cite[\S 3.1]{TypeK-classification} $\Aut(S\times E)=\Aut(S)\times\Aut(E).$ Hence every $g\in G$ is a product $(g_S,g_E).$ The projections $G\to G_S\subset\Aut(S)$ and $G\to G_E\subset\Aut(E)$ are injective. 
\end{proposition}

\begin{proposition}\label{prop:HK2}
    \cite[Prop. 3.3]{TypeK-classification} Let $H=\ker(G\to\GL(H^{2,0}(S)))$ and $\iota\in G\setminus H$. Then $\iota$ has order two and $G=H\rtimes \langle\iota\rangle$ with $\iota h\iota=h^{-1}$ for all $h\in H$. Every $g\in G\setminus H$ acts on $S$ by an Enriques involution. After a choice of origin on $E$, $\iota_E=-1_E$ and $H_E$ is a group of translations.
\end{proposition} 

\begin{proposition}[{\cite[Lemma 3.14, Prop. 4.4, Lemma 4.1]{TypeK-classification}}]\label{prop:HK3}
$H^2(S,\Z)^G=\NS(S)^G$ is as follows:
\begin{center}\footnotesize\setlength{\tabcolsep}{2pt}
\begin{tabular}{lcccccccc}\toprule
$G$ & $C_2$ & $C_2^2$ & $C_2^3$ & $D_6$ & $D_8$ & $D_{10}$ & $D_{12}$ & $C_2\times D_8$\\ \midrule
$\NS(S)^G$ & $U(2)\oplus E_8(-2)$ & $U(2)\oplus D_4(-2)$ & $U(2)\oplus\langle-4\rangle^2$ & $U(2)\oplus A_2(-2)$ & $U(2)\oplus\langle-4\rangle$ & $U(2)$ & $U(2)$ & $U(2)$\\
$\rk$ & 10 & 6 & 4 & 4 & 3 & 2 & 2 & 2\\ \bottomrule
\end{tabular}
\end{center}
The $G$-invariant ample cone is a non-empty open cone in $\NS(S)^G_\R$, and $h^{1,1}(X)=\rk\NS(S)^G+1$.
\end{proposition}

\subsection{Stability conditions on K3 surfaces}\label{subsec:stability-K3}
Let $S$ be a K3 surface. For $F\in\Db(S),$ let $v(F)=\ch(F)\sqrt{\td(S)}=(\rk F,c_1(F),\ch_2(F)+\rk(F))\in\Knum(S)$ be its Mukai vector. The Mukai pairing is defined by $\langle (r,c,s),(r',c',s') \rangle=c\cdot c'-rs'-r's.$ For $\omega,B\in\NS(S)_{\R}$ with $\omega$ ample, set
\begin{equation}\label{Eq: 2.3.1}
    Z_{\omega,B}(F)=\langle e^{B+i\omega},v(F)\rangle=B\cdot c-s-\frac{r}{2}(B^2-\omega^2)+i(\omega\cdot c-rB\cdot\omega),
\end{equation}
Let $\Coh^{\omega,B}(S)=\langle\FF_{\omega,B}[1],\TT_{\omega,B} \rangle_{\mathrm{ext}}$. If $Z_{\omega,B}(\delta)\notin\R_{\leq 0}$ for every $\delta\in\Knum(S)$ with $\langle\delta,\delta\rangle=-2$ and $\rk\delta>0$, then $\sigma_{\omega,B}=(\Coh^{\omega,B}(S),Z_{\omega,B})$ is a full numerical stability condition, and every skyscraper sheaf $\OO_x$ is $\sigma_{\omega,B}$-stable of phase one with $Z_{\omega,B}(\OO_x)=-1$ \cite{Bri-K3}. It is rational when $\omega,B\in\NS(S)_\Q$. We write $U(S)\subset\NS(S)_{\Q}\times\NS(S)_{\Q}$ for the set of pairs $(\omega,B)$ with $\omega$ ample satisfying this condition.

\begin{lemma}\label{lem: U(S)}
    \cite{Bri-K3} $\{(\omega,B)\ |\ \omega^2>2\}\subseteq U(S)$ with equality when $B=0.$.
\end{lemma}

\subsection{Liu's construction on products}
Let $Y$ be a smooth projective variety with a rational numerical stability condition $\sigma=(\AA,Z),$ i.e. image of $Z$ is in $\Q\oplus\Q i$, with respect to $(\Lambda,v).$ By \cite[Lemma 12.12]{Families-paper} $\AA$ is Noetherian. Let $E$ be an elliptic curve, $x_0\in E$, $L=\OO_E(x_0)$ on $E,$ and $p:Y\times E\to Y$, $q:Y\times E\to E$ be the projections. For $x\in E$ let $i_x:Y\to Y\times E,$ $y\mapsto (y,x).$ 
\begin{lemma}\label{Lem: central-charge-first-one}
    For $F\in\Db(Y\times E)$ and all $n\in \Z,$
    \[
    Z(p_*(F\otimes q^*L^n))=Z(p_*F)+nZ(i_x^*F),\text{ for every $x\in E.$}
    \]
    The classes $v(p_*F)$ and $v(i_x^*F)$ only depend on the numerical class of $F$ and $v(i_x^*F)$ is independent of $x.$
\end{lemma}
\begin{proof}
    Tensor the short exact sequence $0\to L^n\to L^{n+1}\to L^{n+1}|_{x_0}\cong\OO_{x_0}\to 0$ with $F$ and apply $p_*.$ We get $p_*(F\otimes q^*L^{n+1})$ as an extension of $p_*(F\otimes\OO_{Y\times\{x_0\}})$ by $p_*(F\otimes q^*L^n).$ Notice, $p_*(F\otimes\OO_{Y\times\{x_0\}})=p_*i_{x_0,*}i_{x_0}^*F=i_{x_0}^*F.$ This completes the proof for $x=x_0.$ Since $p_*$ and $i_x^*$ have adjoints, they preserve numercal equivalence. So $v(p_*F)$ and $v(i_x^*F)$ only depend on the numerical class of $F$. 

    Fix $A\in\Db(Y)$ and set $G=F\otimes p^*A^\vee$. Since $p\circ i_x=\mathrm{id}_Y,$ we get $i_x^*G\cong A^\vee\otimes i_x^*F$. Hence,
    \[
    \chi(A,i_x^*F)=\chi(Y,A^\vee\otimes i_x^*F)=\chi(Y,i_x^*G).
    \]
    Applying derived cohomology and base change to $q$ (which is flat and proper),
    \[
    q_*G\otimes^{\L}k(x)\cong R\Gamma(Y\times\{x\},G|_{Y\times\{x\}})=R\Gamma(Y,i_x^*G).
    \]
    This implies $\chi(A,i_x^*F)$ is the Euler characteristic of the derived fiber at $x$ of $q_*G$. Since $E$ is connected, by semicontinuity, it is constant. This shows $\chi(A,i_x^*F)$ is independent of $x$ for every $A.$ Since a class in $\Knum(Y)$ is determined by its pairings with all $A,$ this proves, $v(i_{x_0}^*F)=v(i_{x}^*F)$ for any $x\in E.$ Since $Z$ factors through $\Knum,$ $Z(i_{x_0}^*F)=Z(i_{x}^*F).$ This completes the proof.
\end{proof}

Now write,
\begin{equation}\label{eq:2.4.1}
    Z(i_x^*F):=a(F)+ic(F),\quad Z(p_*F):=b(F)+id(F),
\end{equation}
where $a,b,c,d:K_0(Y\times E)\to\Q$ are group homomorphisms. Thus, $Z(p_*(F\otimes q^*L^n))=(a(F)n+b(F))+i(c(F)n+d(F)).$

We now follow the construction of stability conditions on $Y\times E$ as in \cite[\S 4]{Liu-product}. 
\[
\AA_E=\{F\in\Db(Y\times E);\ p_*(F\otimes q^*L^n)\in\AA\text{ for }n\gg 0\}
\]
is the heart of a bounded Noetherian $t$-structure on $\Db(Y\times E)$. It is independent of the choice of the degree one line bundle $L$ on $E.$ 

For $t\in\Q_{t>0}$, $Z_t(F)=a(F)t-d(F)+ic(F)t$ is a weak stability function. The pair $(\AA_E,Z_t)$ admits HN property and one constructs the heart $\AA_{E}^t=\langle \FF_t[1],\TT_t\rangle_{\mathrm{ext}},$ where $\FF_t\subset\AA_E$ consists of all the objects whose non-zero subobjects $F'$ satisfy $c(F')\neq 0$ and $a(F')t-d(F')\geq 0$, and $\TT_t$ of all objects whose non-zero quotients $F''$ satisfy $c(F'')=0$ or $a(F'')t-d(F'')<0.$ For $s,t\in\Q_{>0}$,
\begin{equation}\label{eq:2.4.2}
    Z^{L,s,t}(F)=c(F)s+b(F)+i(d(F)-a(F)t),\ \sigma^{L,s,t}=(\AA_E^{t},Z^{L,s,t})
\end{equation}
is a rational numerical stability condition on $\Db(S\times E)$ \cite[Theorem 4.7, 5.9]{Liu-product}, with respect to 
\[
\Lambda_E=\Lambda\oplus\Lambda/{\ker Z},\ v_E(F)=\left(v(i_x^*F),[v(p_*F)]\right)
\]
\cite[Lemma 5.2,5.7, Remark 5.8]{Liu-product}. The lattice $\Lambda_E$ is the quotient of the lattice $\Lambda\oplus\Lambda,$ $F\mapsto (v(i_x^*F),v(p_*F)).$ The $a,b,c,d$ depend on the choice of the line bundle $L$.

\subsection{The central charge on $S\times E$}
Now let $Y=S$ be a K3 surface. Let $(\omega,B)\in U(S),$ and $\AA_S=\Coh^{\omega,B}(S),Z_S=Z_{\omega,B},$ where $Z_S(\gamma)=\langle e^{B+i\omega},\sqrt{\td(S)}\gamma\rangle.$ Let,
\[
\sigma_{\omega,B,s,t}:=(\sigma_{\omega,B})^{L,s,t},\ \Lambda_E=\Knum(S)\oplus\Knum(S)/\ker Z_{\omega,B}.
\]

\begin{lemma}\label{lem:Z-L,s,t}
    For $F\in\Db(S\times E)$, let $\alpha=\ch(i_x^*F)$ and $\beta=\ch(p_*F)$ in $H^*(S,\Q)$, and write $\alpha=\alpha_0+\alpha_1+\alpha_2$ and $\beta=\beta_0+\beta_1+\beta_2$ for the degree $0,2,4$ terms.  Then,
    \begin{align*}
        \Re Z^{L,s,t}(F)&=-\beta_2-\beta_0\left(1+\frac{B^2-\omega^2}{2}\right)+B\cdot\beta_1+s\omega\cdot\alpha_1-s\alpha_0B\cdot\omega,\\
        \Im Z^{L,s,t}(F)&= \omega\cdot\beta_1-\beta_0B\cdot\omega+t\alpha_2+t\alpha_0\left(1+\frac{B^2-\omega^2}{2}\right)-tB\cdot\alpha_1.
    \end{align*}
\end{lemma}
\begin{proof}
      By Lemma \ref{Lem: central-charge-first-one}, $Z(p_*(F\otimes q^*L^n))=nZ(i_x^*F)+Z(p_*F)=(an+b)+i(cn+d).$ Comparing with Eq. \ref{Eq: 2.3.1}, we get,
\begin{align*}
    a&=B\cdot\alpha_1-\alpha_2-\alpha_0\left(1+\frac{B^2-\omega^2}{2}\right),\ c=\omega\cdot\alpha_1-\alpha_0B\cdot\omega,\\
    b&=B\cdot\beta_1-\beta_2-\beta_0\left(1+\frac{B^2-\omega^2}{2}\right),\ d=\omega\cdot\beta_1-\beta_0B\cdot\omega.
\end{align*}
By \eqref{eq:2.4.2}, $\Re Z_E^{s,t}(F)=b+sc$ and $\Im Z_E^{s,t}(F)=d-ta.$ Substituiting the values of $a,b,c,d$ completes the proof.
\end{proof}

By K\"unneth formula, 
\[
H^n(S\times E,\Q)\cong\bigoplus_{i+j=n}H^i(S,\Q) \otimes H^j(E,\Q).
\]
Hence, $H^{2k}(S\times E,\Q)\cong H^{2k}(S,\Q)\otimes H^0(E,\Q)\oplus H^{2k-1}(S,\Q)\otimes H^1(E,\Q)\oplus H^{2k-2}(S,\Q)\otimes H^2(E,\Q).$ $S$ has no odd cohomology, which implies that for every $F\in\Db(S\times E)$, we have,
\begin{equation}\label{eq:2.5.1}
    \ch(F)=p^*\alpha+p^*\beta\cdot f,
\end{equation}
where $\alpha,\beta\in H^{*}(S,\Q)$, and $f=q^*([\mathrm{pt}]).$ Moreover, $\alpha=\ch(i_x^*F)$ and $\beta=p_*\ch(F).$ Since $\td(E)=1,$ by Grothendieck-Riemann-Roch for $p$, $p_*\ch(F)=\ch(p_*F).$

\subsection{Inducing stability conditions}\label{sec:inducing}
Let $G$ be a finite group acting freely on a smooth projective variety $Y'$. Let $\Lambda$ be a finite rank free abelian group equipped with a $G$-action, and let $v:K_0(X)\to\Lambda$ be a $G$-equivariant homomorphism. Let $\pi:Y'\to [Y'/G]$ be the canonical morphism to the quotient stack. For a $G$-invariant $\sigma=(\AA,Z)\in\Stab_{\Lambda,v}(Y')$, set
\[
\AA^G=\{F\in\Db([Y'/G])\ |\ \pi^*F\in\AA\},\ Z^G=Z\circ\pi^*,\ \sigma^G=(\AA^G,Z^G).
\]

\begin{theorem}\label{thm:MMS}
    \cite[Theorem 2.16]{LMPSZ}, \cite{MMS}, \cite[Theorem 4.3]{Perr-Shah} $\sigma^G$ is a stability condition on $\Db([Y'/G])$ with respect to $(\Lambda,v\circ\pi^*),$ and also with respect to the invariant lattice $(\Lambda^G,v^G)$ through which $v\circ\pi^*$ factors. If $\sigma$ is numerical, then so is $\sigma^G.$ An object $F$ is $\sigma^G$-semistable if and only if $\pi^*F$ is $\sigma$-semistable, and in that case, they share the same phase.
\end{theorem}

If $G$ acts freely, then $Y'/G:=X$ is smooth projective and $\Db(X)\cong\Db([Y'/G])$ and the pullback $\pi^*$ along the \'etale cover $\pi:Y'\to X$ is exact. For $F\in\Db(X)$, one has $\pi^*\pi_*F\cong\bigoplus_{g\in G}g^*F.$

\section{Stability conditions on CY3-folds of type K}

Suppose we are in the setup of the previous section. We closely follow the ideas in \cite[\S 4]{Perr-Shah}. In particular, look at \cite[Lemma 4.8, 4.9]{Perr-Shah}. 
\begin{lemma}
    Let $\gamma=(g,h)\in\Aut(Y)\times\Aut(E)$ act on $Y\times E$, and let $\sigma=(\AA,Z)\in\Stab(Y)$ be $g$-invariant. Then $\gamma^*\AA_E^t=\AA_E^t$ and $Z^{L,s,t}\circ \gamma=Z^{L,s,t}$ for all $s,t>0.$
\end{lemma}
\begin{proof}
    Notice, $\gamma$ is an automorphism such that $p\gamma=gp$ and $q\gamma=hq$. Thus, $p_*\gamma^*=g^*p_*$ and $\gamma^*q^*=q^*h^*$, so,
    \[
    p_*(\gamma^*F\otimes q^*L^n)=p_*\gamma^*(F\otimes q^*(h_*L)^n)=g^*p_*(F\otimes q^*(h_*L)^n).
    \]
    Since $\sigma$ is $g$-invariant $g^*\AA=\AA.$ Therefore, $\gamma^*F\in\AA_E$ if and only if $F$ lies in the heart defined with $h_*L$ in place of $L$. $\AA_E$ is independent of the choice of degree one line bundle. This proves $\gamma^*\AA_E=\AA_E.$ Next, $p_*\gamma^*F=g^*p_*F$ and $i_x^*\gamma^*F=g^*i_{hx}^*F.$ Thus, by $Z\circ g^*=Z$ and Lemma \ref{Lem: central-charge-first-one}, the functions $a,b,c,d$ are $\gamma$-invariant. Hence, the central charge, the torsion pair, and the tilt are fixed by $\gamma.$ This completes the proof.
\end{proof}

\begin{lemma}\label{lemma:grp-K3}
    Let $G$ act on a K3 surface $S$ through $G_S\subset \Aut(S).$ For $(\omega,B)\in U(S)$ with $\omega,B\in\NS(S)_\Q^G$, $\sigma_{\omega,B}$ is $G_S$-invariant. Such $(\omega,B)$ exist for all the eight cases.
\end{lemma}
\begin{proof}
    The first statement is immediate from the tilting construction. By Proposition \ref{prop:HK3}, there exists $G$-invariant ample classes $\omega,$ and $(\omega,B)\in U(S)$ whenever $\omega^2>2$ by Lemma \ref{lem: U(S)}.
\end{proof}

\begin{proposition}\label{mainprop-Liu-G-invariant}
    Let $X$ be a type K threefold. For $(\omega,B)\in U(S)$ with $\omega,B\in\NS(S)_{\Q}^G$, and $s,t\in\Q_{>0},$ the stability condition $\sigma_{\omega,B,s,t}$ on $\Db(S\times E)$ is $G$-invariant. Moreover, $v_E:K_0(S\times E)\to\Lambda_E$ is $G$-equivariant for the action of $G$ on $\Lambda_E$ induced by $g_S^*$ on $\Knum(S).$
\end{proposition}
\begin{proof}
    Every $g\in G$ is of the form $(g_S,g_E)$ (Prop. \ref{prop:HK1}). Then $G$-invariance of $\sigma_{\omega,B,s,t}$ follows from the the previous two lemmas. Notice that the action on $\Lambda_E$ is well-defined, because $g_S^*$ preserves $\ker Z_{\omega,B}.$ Lastly, $i_x^*g^*=g_S^*i_{g_E x}^*$ and $p_*g^*=g_S^*p_*,$ which by Lemma \ref{Lem: central-charge-first-one} implies $v_E(g^*F)=g^*v_E(F).$ Hence $v_E$ is $G$-equivariant. 
\end{proof}

\begin{theorem}\label{mainthm}
    In the setting of the previous proposition,
    \[
    \sigma_{\omega,B,s,t}^G=\left(\{F\in\Db(X);\ \pi^*F\in\AA_E^t\},Z^{L,s,t}\circ\pi^*\right)
    \]
    is a numerical Bridgeland stability condition with respect to $(\Lambda_E,v_E\circ\pi^*).$ An object $F\in\Db(X)$ is $\sigma^G$-semistable if and only if $\pi^*F$ is $\sigma_{\omega,B,s,t}$-semistable. The family has dimension $2\rk\NS(S)^G+2=2h^{1,1}(X).$
\end{theorem}
\begin{proof}
    The proof of the first two statements follows from Prop. \ref{mainprop-Liu-G-invariant} and Theorem \ref{thm:MMS}.

    For $F\in\Db(X),$ $\ch(\pi^*F)=p^*\alpha+p^*\beta\cdot f,$ where $\alpha=\alpha_0+\alpha_1+\alpha_2,\beta=\beta_0+\beta_1+\beta_2\in H^*(S,\Q)$. Set $Z^G(F)=(Z^{\omega,B,s,t})^G(F)=Z^{\omega,B,s,t}(\pi^*F).$ For any $g=(g_S,g_E)\in G$, we get,
    \[
    g^*(p^*\alpha+p^*\beta\cdot f)=p^*g_S^*\alpha+p^*g_S^*\beta\cdot f.
    \]
    Since $g^*\pi^*=\pi^*$, $g_S^*\alpha=\alpha,g_S^*\beta=\beta,$ so $(\alpha,\beta)\in V\oplus V,$ where $V=\Q\oplus\NS(S)_{\Q}^G\oplus\Q\subset H^*(S,\Q)^G.$ 
    
    Conversely, let $A\in\Db(S\times E)$ with $\ch(A)=p^*\alpha'+p^*\beta'\cdot f$ and set $F=\pi_*A.$ We have,  $\pi^*F=\pi^*\pi_*A\cong\bigoplus_{g\in G}g^*A,$ so $\ch(\pi^*F)=p^*(\sum_g g_S^*\alpha')+p^*(\sum_g g_S^*\beta')\cdot f$. The map $\gamma\mapsto\sum_{g}g_S^*\gamma$ maps $\Q\oplus\NS(S)_\Q\oplus\Q$ onto $V.$ Take $A=p^*A_1,$ which gives $(\alpha',\beta')=(\ch A_1,0)$ and if we take $A=i_{e*}A_2,$ then $(\alpha',\beta')=(0,\ch A_2).$ Since the Chern characters of objects on $\Db(S)$ span $\Q\oplus\NS(S)_\Q\oplus\Q,$ the pairs $(\alpha,\beta)$ span $V\oplus V.$

    By Lemma \ref{lem:Z-L,s,t}, $Z^G$ depends on $F$ by $\alpha,\beta$ and the dependence is linear. In other words, there's a linear map
    \[
    l_{\omega,B,s,t}:V\oplus V\to \C,\ Z^G(F)=l_{\omega,B,s,t}(\alpha,\beta).
    \]
    Suppose $(Z_{\omega,B,s,t})^G=(Z_{\omega',B',s',t'})^G$ as functions on $K_0(X).$ So $l_{\omega,B,s,t}=l_{\omega',B',s',t'}$ on $V\oplus V$. Now we compare the linear map on each summand of $V\oplus V$.
    \begin{enumerate}
        \item On the $\beta_1$-summand: $\Re l=B\cdot (-)$ and $\Im l=\omega\cdot (-).$ So $B\cdot D=B'\cdot D$ for all $D\in\NS(S)^G_\Q.$ Hence $B=B'$ by the non-degeneracy of the intersection form on $\NS(S)^G.$ Similarly, $\omega=\omega'.$
        \item On the $\alpha_2$-summand: $\Im l=t\cdot (-)\implies t=t'.$
        \item On the $\alpha_1$-summand: $\Re l=s\omega\cdot (-),$ so $s\omega\cdot D=s'\omega\cdot D$ for all $D.$ Take $D=\omega.$ Since $\omega^2>0,$ we get $s=s'.$
    \end{enumerate}
    Hence $(\omega,B,s,t)=(\omega',B',s',t')$.
    
    The parameters range over $\{(\omega,B)\in U(S)\cap (\NS(S)^G_{\Q})^2\}\times\Q_{>0}^2,$ which is dense in an open subset of $(\NS(S)^G_{\R})^2\times\R^2_{>0}$. This proves the dimension count and the equality $2\rk\NS(S)^G+2=2h^{1,1}(X)$ follows from Prop. \ref{prop:HK3}.
\end{proof}

\section{Objects supported on a fiber}\label{sec:fiber}
Let $Y$ be smooth projective and $\sigma=(\AA,Z)$ a rational numerical stability condition on $Y$ and $\sigma'=\sigma^{L,s,t}=(\AA',Z')\in\Stab(Y\times E).$ We follow the ideas from \cite[\S 3]{LMPSZ}. Fix $e\in E$ and let $\DD_e=\Db_{Y\times\{e\}}(Y\times E)$ be the full subcategory of objects set-theoretically supported on $Y\times \{e\}.$

\begin{lemma}\label{lem1:fiber}
    \begin{enumerate}
        \item For any $K\in\DD_e$ and every line bundle $M$ on $E,$ we have $K\otimes q^*M\cong K$.
        \item If $K\in\DD_e$ and $p_*K=0$, then $K=0.$
        \item For any $F\in\Db(Y)$, $i_e^*i_{e*}F\cong F\oplus F[1]$. This implies $\Hom(i_{e*}F,i_{e*}F)\cong\Hom(F,F)\oplus\Hom(F[1],F).$
    \end{enumerate}
\end{lemma}
\begin{proof}
    \begin{enumerate}
        \item Let $U\subset E$ be an open neighborhood of $e$, such that $M|_{U}$ is trivial. Let $j:Y\times U\hookrightarrow Y\times E$ be the inclusion. The restriction to $Y\times U$ induces an equivalence $\DD_e\cong\Db_{Y\times \{e\}}(Y\times U).$ Since $M|_U$ is trivial, $j^*(K\otimes q^*M)\cong j^*K.$ Thus $K\otimes q^*M\cong K.$
        
        \item Since $p$ has relative dimension one, $\R^ip_*=0$ for $i\neq 0,1,$ so $p_*$ has cohomological amplitude $[0,1].$ Let $m$ be minimal with $\mathcal{H}^m(K)\neq 0$. Apply $p_*$ to the triangle $\mathcal{H}^m(K)[-m]\to K\to\tau_{\geq m+1}K.$ After pushforward the third term lies in $\Db(Y)^{\geq m+1}$ and the first term lies in $\Db(Y)^{\geq m}$. Thus the long exact sequence gives,
        \[
        \mathcal{H}^m(p_*K)\cong \R^0p_*\mathcal{H}^m(K).
        \]
        The cohomology sheaf $\mathcal{H}^m(K)$ is supported in $Y\times \{e\},$ so it is the pushforward of a non-zero sheaf $G$ on $Y\times T$ for some finite subscheme $T\subset E$ supported at $e.$ So $\R^0p_*\mathcal{H}^m(K)=(p|_{Y\times T})_*G.$ The morphism $p|_{Y\times T}$ is affine, so the pushforward is non-zero, since $G$ is non-zero. Thus $\mathcal{H}^m(p_*K)\neq 0$, which implies $p_*K\neq 0$.
        
        \item Since $p\circ i_e=\mathrm{id}_Y$, projection formula gives $i_{e*}F\cong p^*F\otimes\OO_{Y\times\{e\}}\cong p^*F\otimes\OO_{Y\times\{e\}}.$ Hence $i_e^*i_{e*}F\cong F\otimes i_e^*\OO_{Y\times\{e\}}.$ Now use the locally free resolution $q^*\OO_E(-e)\xrightarrow{s}\OO_{Y\times E}.$ Pulling back by $i_e$ gives $\OO_Y\to\OO_Y$ in degrees $-1,0.$ The differential is $s|_{Y\times\{e\}}=0.$ Hence $i_e^*\OO_{Y\times\{e\}}\cong\OO_Y\oplus\OO_Y[1]$. Thus $i_e^*i_{e*}F\cong F\oplus F[1].$ The last statement follows from the adjunction $\Hom(i_{e*}F,i_{e*}F)\cong\Hom(i_e^*i_{e*}F,F).$
    \end{enumerate}
\end{proof}

\begin{lemma}\label{lem:finite-supp}
    Let $K\in\Db(Y\times E)$ such that $K\otimes  q^*M\cong K$ for all $M\in\Pic^0(E).$ Then $q(\Supp K)$ is finite and $K=\bigoplus_{x\in q(\Supp K)}K_x$, where $K_x$ is supported on $Y\times \{x\}.$
\end{lemma}
\begin{proof}
    Follows from \cite[Lemma 3.3]{LMPSZ}.
\end{proof}

\begin{lemma}\label{lem2-fiber}
    $\sigma'$ is invariant under tensoring with $q^*M,$ for $M\in \Pic^0(E).$ If $K\in\DD_e$, then every Harder--Narasimhan factor of $K$ with respect to $\sigma'$ lies in $\DD_e.$ Hence, $\sigma_e':=(Z',\PP(\phi)\cap\DD_e)$ is a stability condition on $\DD_e$ with respect to $(\Lambda_E,v_E|_{\DD_e}).$ Moreover, $K\in\DD_e$ is $\sigma'$-semistable if and only if it is $\sigma_e'$-semistable.
\end{lemma}
\begin{proof}
    Following \cite[Theorem 2.23]{LMPSZ}, where we take $Y=X=S\times E$, $\AA'$ is invariant under the $\Pic^0(E)$ action (it is Noetherian because $\sigma'$ is rational). Since $v_{\sigma'}$ factors through $\Knum(S\times E)$ (Lemma \ref{Lem: central-charge-first-one}), \cite[Theorem 2.23 (2)]{LMPSZ} is satisfied by \cite[Remark 2.24]{LMPSZ}. This proves $\sigma'$ is invariant under the $\Pic^0(E)$-action. Rest follows from \cite[Lemma 3.4, Prop. 3.1]{LMPSZ} and Lemma \ref{lem1:fiber}.
\end{proof}

\begin{lemma}\label{lem3-fiber}
    Let $K\in\DD_e$. Then,
    \begin{enumerate}
        \item $a(K)=c(K)=0$ and $Z'(K)=Z(p_*K).$ $K\in\AA_E$ if and only if $p_*K\in\AA.$
        \item $\AA'\cap\DD_e=\AA_E\cap\DD_e=:\AA_e=\{K\in\DD_e;\ p_*K\in\AA\}.$ This is an abelian subcategory of $\AA_E.$ Moreover, it is closed under subobjects, quotients, and extensions in $\AA_E.$
        \item $p_*:\AA_e\to\AA$ is exact and faithful. $i_{e*}:\AA\to\AA_e$ is exact and both functors preserve central charges.
    \end{enumerate}
    In particular, the restricted stability condition $\sigma_e'$ from the previous lemma, is $(Z\circ p_*,\AA_e)$, and is independent of $s,t,L.$
\end{lemma}
\begin{proof}
    \begin{enumerate}
        \item The first two are clear. Recall $\AA_E=\{F\in\Db(S\times E);\ p_*(F\otimes q^*L^n)\in\AA\}.$ For $K\in\DD_e,$ $p_*(K\otimes q^*L^n)=p_*K$ (Lemma \ref{lem1:fiber} (1)).
        \item Let $U=E\setminus\{e\}$ and $j:Y\times U\hookrightarrow Y\times E.$ Since $\AA_E$ is $E$-local, the pullback $j^*$ is $t$-exact for the hearts $\AA_E$ and $\AA_U.$ Thus $j^*H_{\AA_E}^i(K)=H_{\AA_U}^i(j^*K)$. If $K\in\DD_e$, $j^*K=0$, hence $j^*H_{\AA_E}^i(K)=0$ for ann $i.$ An object lies in $\DD_e$ when the restriction to $Y\times U$ vanishes. This implies $H_{\AA_E}^i(K)\in\DD_e.$
        
        Take $K,K'\in\AA_e.$ For $\phi:K\to K'$ in $\AA_E,$ the cone $\mathcal{C}$ lies in $\DD_e$ and by above, $\ker\phi=H_{\AA_E}^{-1}(\mathcal{C})$ and $\mathrm{coker}\phi=H_{\AA_E}^0(\mathcal{C})$ both lie in $\DD_e.$ Hence, $\AA_e\subset\AA_E$ is a full subcategory closed under kernels and cokernels, hence abelian. An extension in $\AA_E$ is a cone of a map between objects in $\DD_e,$ so it lies in $\DD_e.$ If $K'\hookrightarrow K$ is a monomorphism in $\AA_E$ with $K\in\DD_e,$ then $j^*K'\hookrightarrow j^*K=0$ is a monomorphism in $\AA_U$ ($t$-exact). Thus, $j^*K'=0$ and so $K'\in\DD_e.$ Similarly for quotients.
        
        By (1), a non-zero object of $\AA_E\cap\DD_E$ has $c=0$. So it cannot be in $\FF_t$ and all its quotients have $c=0$, so it lies in $\TT_t.$ This implies $\FF_t\cap\DD_e=0$ and since $\TT_t\subset\AA',$ $\AA_E\cap\DD_e\subset\AA'\cap\DD_e.$  Conversely, let $K\in\AA'\cap\DD_e=\langle\FF_t[1],\TT_t\rangle.$ We have $H^{-1}_{\AA_E}(K)\in\FF_t,H^{0}_{\AA_E}(K)\in\TT_t.$ Again, by the first argument above, $H^{-1}_{\AA_E}(K)\in\DD_e,$ so $H^{-1}_{\AA_E}(K)\in\DD_e\cap\FF_t=0.$ Therefore, $K=H^{0}_{\AA_E}(K)\in\AA_E.$ So $K\in\AA_E\cap\DD_e.$
        
        \item Exactness is immediate. Faithfulness follows from Lemma \ref{lem1:fiber} (2). They preserve central charge by $(1).$
    \end{enumerate}
\end{proof}

Let $K\in\AA_e$ of phase $\phi.$ By Lemma \ref{lem2-fiber} and Lemma \ref{lem3-fiber}, $K$ is $\sigma'$-semistable if and only if every subobject in $\AA_e$ has phase $\leq \phi.$

\begin{theorem}\label{mainthm: stability-fiber}
    Let $F\in\Db(Y).$ Then $F$ is $\sigma$-semistable (resp. stable) if and only if $i_{e*}F$ is $\sigma'$-semistable (resp. stable). They share the same phase.
\end{theorem}
\begin{proof}
    After shifting, we may assume $F\in\AA$. They share the same phase by Lemma \ref{lem1:fiber} (3). If $F'\subset F$ destabilizes $F$, then $i_{e*}F'\subset i_{e*}F$ destabilizes $i_{e*}F$ in $\AA_e.$ Conversely, let $0\neq K\subset i_{e*}F$ is a subobject in $\AA_e.$ Then $0\neq p_*K\subset F$ in $\AA$ (by Lemma \ref{lem1:fiber} (2)), and $p_*$ preserves central charges. If $F$ is $\sigma$-semistable, $\phi(K)\leq\phi(F)$ for every subobject of $i_{e*}F.$  Thus $i_{e*}F$ is $\sigma'$-semistable by the remark before the theorem.

    Let $i_{e*}F$ be $\sigma'$-stable. If $F$ is semistable but not stable, a proper non-zero subobject $F'\subset F$ of phase $\phi$ gives a proper non-zero subobject $i_{e*}F'\subset i_{e*}F$ in $\PP(\phi),$ contradictory to the stability of $i_{e*}F.$

    \emph{(subobjects in $\DD_e$)} Now let $F$ be stable and $K\subset i_{e*}F$ be an $\AA_e$-subobject of phase $\phi.$ Then $p_*K\subset F$ has phase $\phi.$ Therefore, $p_*K=0$ or $p_*K=F$. If $p_*K=0$ then $K=0$ by Lemma \ref{lem1:fiber} (2). Likewise, if $p_*K=F$, then $p_*(i_{e*}F/K)=0$, so $K=i_{e*}F.$

    \emph{(general subobjects)} Suppose $i_{e*}F$ is not $\sigma'$-stable and let $K\subsetneq i_{e*}F$ be a subobject in $\PP(\phi).$ For any $M\in\Pic^0(E)$, $i_{e*}F\otimes q^*M\cong i_{e*}F$ and $\PP(\phi)\otimes q^*M=\PP(\phi)$ (Lemma \ref{lem:finite-supp}). Let $K_M:=\operatorname{im}(K\otimes q^*M\to i_{e*}F)$ be the subobject in $\PP(\phi)$. Since $\PP(\phi)$ has finite length, a sum of simple subobjects of a finite-length object is a finite direct sum of some of them. This implies that the sum $K':=\sum_M K_M=\bigoplus_{j=1}^r K_{M_j}.$ Notice, $K_M\otimes q^*N\cong K_{M\otimes N}.$ Thus, $K'\otimes q^*N\cong K'$ for any $N\in\Pic^0(E).$ By Lemma \ref{lem2-fiber}, $K'=\bigoplus_x K'_x$, where $K_x'$ is supported on $Y\times\{x\}.$ Since $\PP(\phi)$ is closed under direct summands, $K_x'\in\PP(\phi).$ Notice, $K_x$ maps to zero in $i_{e*}F$, for $x\neq e,$ since they have disjoint support. But $K_x\hookrightarrow K'\hookrightarrow i_{e*}F$ would give a monomorphism, which forces $K'_x\hookrightarrow K'$ to be zero for any $x\neq e.$ So, $K'\in\DD_e.$ Since the cokernel of $K'\to i_{e*}F$ also lies in $\DD_e,$ we conclude that $K'$ is a subobject of $i_{e*}F$ in $\AA_e$ of phase $\phi.$ By the previous step, $K'=i_{e*}F$. If $r=1,$ then $K'=K_{M_1}=i_{e*}F$ is a simple object in $\PP(\phi).$ This is contradictory to the proper non-zero subobject $K\subsetneq i_{e*}F.$ If $r\geq 2,$ then $i_{e*}F$ is decomposable. Since $F$ is stable, $\operatorname{End}(F)=\C$ and $\Hom(F[1],F)=0$, because $F\in\AA.$ By Lemma \ref{lem1:fiber} (3), we get $\operatorname{End}(i_{e*}F)=\C.$ Hence decomposability of $i_{e*}F$ again leads to a contradiction.
\end{proof}

\subsection{Consequences on $X$}
Let $X=(S\times E)/G$ be of type $K.$ Then $\overline{q}:X\to\P^1=E/G_E$ is the K3 fibration. Let $E^\circ\subset E$ be the open subset on which $G_E\cong G$ acts freely. For $e\in E^\circ$, the fiber $X_{e}$ is smooth. Notice, $\pi^{-1}(X_e)=\bigsqcup_{g_E\in G_E}S\times \{g_Ee\}$ and $\pi\circ i_e:S\to X_e$ is an isomorphism. 

Set $j_e:S\cong X_e\hookrightarrow X.$ We have $g^{-1}\circ i_e=i_{g_E^{-1}e}\circ g_S^{-1}.$   Recall the identity from \S \ref{sec:inducing},
\[
\pi^*\pi_*F\cong\bigoplus_{g\in G}g^*F.
\]

\begin{lemma}\label{lem:decomp}
    $\pi^*j_{e*}F\cong\bigoplus i_{g_Ee*}((g_S^{-1})^*F).$
\end{lemma}
\begin{proof}
    $j_e=\pi\circ i_e.$ So,
    \begin{align*}
        \pi^*j_{e*}F&\cong  \pi^*\pi_*i_{e*}F
        \cong \bigoplus_{g\in G}g^*i_{e*}F
        \cong \bigoplus_{g\in G}(g^{-1}\circ i_{e})_*F
        \cong \bigoplus_{g\in G}(i_{g_E^{-1}e}\circ g_S^{-1})_*F
        \cong \bigoplus_{g\in G}i_{g_E^{-1}e*}g_{S}^*F\\
        &\cong \bigoplus i_{g_Ee*}((g_{S}^{-1})^*F).
    \end{align*}
\end{proof}

\begin{theorem}\label{cor:main1}
    Let $e\in E^\circ$ and $F\in\Db(S).$ Let $\omega,B\in\NS(S)^G.$ Then $j_{e*}F$ is $\sigma^G_{\omega,B,s,t}$-semistable (resp. stable) if and only if $F$ is $\sigma_{\omega,B}$-semistable (resp. stable). We also have $Z^G(j_{e*}F)=|G|Z_{\omega,B}(F).$
\end{theorem}
\begin{proof}
    $\omega,B\in\NS(S)^G$ implies $\sigma_{\omega,B}$ is $G_S$-invariant (Lemma \ref{lemma:grp-K3}). So $(g_S^{-1})^*F$ is $\sigma_{\omega,B}$-(semi)stable if and only if $F$ is, and they share the same phase. If $F$ is $\sigma_{\omega,B}$-(semi)stable, then by Theorem \ref{mainthm: stability-fiber}, the right side of Lemma \ref{lem:decomp} is a direct sum of (semi)stable objects of the same phase. If $F$ is not semistable, then the summand for $g=1$ is not semistable. Since a direct sum is semistable if and only if the summands are semistable of same phase, the proof of equivalence of semistability follows from Theorem \ref{mainthm}.

    Write the decomposition in $\pi^*j_{e*}F\cong  \bigoplus_{g\in G}A_g$, where $g=(g_E,g_S)$. Let $F$ is stable. Then the summands are stable of same phase and pairwise non-isomorphic, since they have distinct supports. The goal is to prove $j_{e*}F$ is stable, i.e. has no proper subobject in $\PP^G(\phi).$ Let $K\subset j_{e*}F$ in $\PP^G(\phi)$. This pulls back to a $G$-invariant subobject of a direct sum of non-isomorphic simples. Hence $\pi^*K=\bigoplus_{g\in I} A_g$, for some $I\subset G.$ The support of $\pi^*K=\bigsqcup_{g\in I} S\times\{g_Ee\}$ is $G$-invariant. But $G$ acts on the fibers transitively, which implies $I=\emptyset$ or $I=G.$ Hence $K=0$ or $K=j_{e*}F$, so $j_{e*}F$ is stable. Now, if $F$ is strictly semistable, there is a non-zero proper subobject $F'\subset F$ of the same phase. We get a short exact sequence $0\to F'\to F\to F/F'\to 0$ in $\PP_S(\phi).$ Now applying the exact functor $j_{e*},$ we get
    \[
    0\to j_{e*}F'\to j_{e*}F\to j_{e*}(F/F')\to 0
    \]
    in $\PP^G(\phi)$ be the previous arguments. From the Lemma \ref{lem:decomp} $\pi^*j_{e*}F'$ contains $i_{e*}F'$ ($g=1$), which is non-zero $(p_*i_{e*}F'=F'\neq 0)$. Thus, $\pi^*j_{e*}F'$ is non-zero. Likewise, the cokernel is non-zero since $F/F'\neq 0.$ This proves $j_{e*}F'$ is a proper non-zero subobject of $j_{e*}F$ and has the same phase. So $j_{e*}F$ is not stable.
\end{proof}

\begin{corollary}\label{cor:main2}
    All skyscraper sheaves $\OO_x$, where $x\in X,$ are $\sigma_{\omega,B,s,t}^G$-stable of phase one, with $Z^G(\OO_x)=-|G|$.
\end{corollary}
\begin{proof}
    Since $G$ acts freely, for any $x\in X,$ $\pi^{-1}(x)=\{g\cdot y\ |\ g\in G\}$ consists of $|G|$ distinct points, where $y=(s,e')\in S\times E.$ Thus, $\pi^*\OO_x\cong\oplus_{g\in G}\OO_{g\cdot y}.$ Each summand has the form $\OO_{(s',e'')}\cong i_{e''*}\OO_{s'}$. Since $(\omega,B)\in U(S)$, $\OO_{s'}$ is $\sigma_{\omega,B}$-stable of phase one with $Z_{\omega,B}(\OO_{s'})=-1$ (Lemma \ref{lem: U(S)}). By Theorem \ref{mainthm: stability-fiber} with $Y=S$, each $\OO_{g\cdot y}$ is $\sigma_{\omega,B,s,t}$-stable of phase one with $Z^{L,s,t}(\OO_{g\cdot y})=-1$. In particular, $\pi^*\OO_x$ is $\sigma_{\omega,B,s,t}$-semistable of phase one. By Theorem \ref{thm:MMS}, $\OO_x$ is $\sigma^G_{\omega,B,s,t}$-semistable of phase one, and $Z^G(\OO_x)=Z^{L,s,t}(\pi^*\OO_x)=-|G|.$

    Now we show that $\OO_x$ is stable. Let $0\neq K\subset\OO_x$ be a subobject in $\PP^G(1)$. Then $\pi^*K\subset\pi^*\OO_x$ in $\PP(1)$. The object $\pi^*\OO_x$ is a direct sum of pairwise non-isomorphic simple objects of $\PP(1)$, since their supports are distinct. Hence $\pi^*K\cong\bigoplus_{z\in I}\OO_z$ for some subset $I\subset\pi^{-1}(x)$, and $\Supp(\pi^*K)=I$. Since $\pi\circ g=\pi$, we have $g^*\pi^*K\cong\pi^*K$ for all $g\in G$, so $I$ is $G$-invariant. As $G$ acts transitively on $\pi^{-1}(x)$, either $I=\emptyset$ or $I=\pi^{-1}(x)$. In the first case $\pi^*K=0$. In the second case, the pullback of the cokernel of $K\to\OO_x$ under $\pi$ is zero. Since $\pi$ is faithful $K=0$ or $K=\OO_x$. Hence $\OO_x$ is $\sigma^G_{\omega,B,s,t}$-stable of phase one.
\end{proof}
\begin{remark}
    Let $w\in\Knum(S)^G.$ For any $F\in\Db(S)$ with $[F]=w,$ and $e\in E^\circ,$ Theorem \ref{cor:main1} identifies the set of $\sigma^G$-stable objects of the form $j_{e*}F$ with the quotient $(M_{S,\sigma_{\omega,B}}(w)\times E^0)/G$, where $g\in G$ acts by $(F,e)\mapsto ((g_S^{-1})^*F,g_Ee).$ Indeed $j_{g_Ee}\circ g_S=j_e$ and $j_{e*}F\cong j_{e'*}F'$ forces $[e]=[e']$ and with Lemma \ref{lem:decomp} implies $(F',e')$ lies in the $G$-orbit of $(F,e).$ 
\end{remark}

\section{Comparison with Li's stability conditions}
Let $X=(S\times E)/G$ be a Calabi--Yau threefold of type K and $\iota:X\hookrightarrow\P^n$ a closed embedding for some $n.$ Let $H=\iota^*\OO(1)$ and $\NN=\NN_{X/\P^n}$. In \cite{Li}, Li constructs stability conditions on $\Db(\P^n)$ with central charge $-\int_{\P^n}e^{-(b+ia)H}\ch$ and shows that for $a\gg 0,$ these restrict to stability conditions on $\Db(X).$ The central charge for the restricted stability condition is \cite[Remark 6.3]{LLLMPSZ:26},
\begin{equation}\label{eq:central-charge:Li}
    Z_{a,b}^{\mathrm{Li}}(F)=-\int_X e^{-(b+ia)H}\td(\NN)^{-1}\ch(F).
\end{equation}

Since $\NS(S\times E)=p^*\NS(S)\oplus q^*\NS(E),$ we can write $\pi^*H=p^*\omega_0+mf$ for $\omega_0$ ample and $m>0.$ Since $c_1(X)=0$, we get $c_1(\NN)=(n+1)H$ and $c_2(\NN)=\binom{n+1}{2}H^2-c_2(X).$ Thus,
\[
\td(\NN)^{-1}=1-\frac{n+1}{2}H+\frac{(n+1)(3n+4)}{24} H^2+\frac{1}{12}c_2(X)+T_3,
\]
where $T_3$ is the degree-$3$ part. For $F\in\Db(X),$ write $\pi^*\ch(F)=p^*\alpha+p^*\beta\cdot f$ with $\alpha=\alpha_0+\alpha_1+\alpha_2$ and $\beta=\beta_0+\beta_1+\beta_2.$ We use these facts for the calculations below: $\pi^*\gamma\cdotp^*\gamma'=\pi^*(\gamma\gamma'),$ $p^*\gamma=0$ if $\deg\gamma>4,$ $|G|\int_X=\int_{S\times E}\pi^*$. Also, $\pi$ is \'etale, so $\pi^*T_X=T_{S\times E}=p^*T_S\oplus q^*T_E.$ $T_E$ is trivial, so $c(S\times E)=p^*c(S)=1+24p^*[\mathrm{pt}]$. This gives $\pi^*c_2(X)=c_2(S\times E)=24p^*[\mathrm{pt}].$
\begin{align*}
    C_1&:=H^2\ch_1=\frac{1}{|G|}(\beta_0\omega_0^2+2m\omega_0\cdot\alpha_1),\quad C_2:=H\ch_2=\frac{1}{|G|}(\omega_0\cdot\beta_1+m\alpha_2),\\
    C_3&:=\ch_3=\frac{\beta_2}{|G|},\quad C_4:=c_2(X)\cdot\ch_1=\frac{24\beta_0}{|G|}\\
    H^3&=\frac{3m\omega_0^2}{|G|},\quad c_2(X)\cdot H=\frac{24m}{|G|}. 
\end{align*}

Let us calculate $C_1$ to showcase the chain of ideas. $(\pi^*H)^2=(p^*\omega_0+mf)^2=p^*\omega_0^2+2mp^*\omega_0\cdot f.$ Now multiplying this by $\pi^*\ch_1=p^*\alpha_1+\beta_0 f,$
\[
(\pi^*H)^2\pi^*\ch_1=p^*\omega_0^2\cdot p^*\alpha_1+\beta_0 p^*\omega_0^2\cdot f+2mp^*(\omega_0\alpha_1)\cdot f+2m\beta_0p^*\omega_0\cdot f^2.
\]
The first term has degree $6$ on $S,$ hence zero. Now integrating,
\begin{align*}
    \int_{S\times E}(\pi^*H)^2\pi^*\ch_1&=\int_{S\times E} \beta_0p^*\omega_0^2\cdot f+\int_{S\times E} 2mp^*(\omega_0\cdot\alpha_1)\cdot f\\
    &=\int_S\beta_0\omega_0^2+\int_S2m\omega_0\cdot\alpha_1\\
    \implies \int_{X}H^2\ch_1&=\frac{1}{|G|}(\beta_0\omega_0^2+2m\omega_0\cdot\alpha_1).
\end{align*}
The rest of them are calculated similarly. Set $\widetilde{\ch}=\td(\NN)^{-1}\ch.$ This means, 
\begin{align*}
    H^2\widetilde{\ch}_1&=C_1-(\tfrac{n+1}{2})H^3r;\ H\widetilde{\ch}_2=C_2-(\tfrac{n+1}{2})C_1+(K_1H^3+\tfrac{1}{12}c_2(X)H)r\\
    \widetilde{\ch}_3&=C_3-(\tfrac{n+1}{2})C_2+K_1C_1+\tfrac{1}{12}C_4+D_3r,
\end{align*}
where $K_1:=\tfrac{(n+1)(3n+4)}{24}, D_3:=\int_X T_3.$ Then one can use this to expand \ref{eq:central-charge:Li},
\begin{align*}
    \Im Z_{a,b}^{\mathrm{Li}}&=aC_2-a(b+\tfrac{n+1}{2})C_1+aK_2r,\quad K_2=(K_1+\tfrac{n+1}{2}b+\tfrac{3b^2-a^2}{6})H^3+\tfrac{1}{12}c_2(X)\cdot H,\\
    \Re Z_{a,b}^{\mathrm{Li}}&=-C_3+(b+\tfrac{n+1}{2})C_2+(-K_1-\tfrac{n+1}{2}b-\tfrac{b^2-a^2}{2})C_1-\tfrac{1}{12}C_4+C_Rr,
\end{align*}
where $C_R$ is a constant depending on $(a,b,n,H).$ The central charge of $\sigma_{\omega,B,s,t}$ is $Z^G(F)=Z^{L,s,t}(\pi^*F),$ and is given by Lemma \ref{lem:decomp}.

\begin{theorem}\label{thm:Li-vs-Liu}
    For all $(\omega,B)\in U(S)$ with $\omega,B\in\NS(S)_{\Q}^G$, all $s,t>0,$ all $a>0$ and $b\in \R,$ $Z^G_{\omega,B,s,t}\neq M\circ Z_{a,b}^{\mathrm{Li}}$ for any $M\in\GL_2^+(\R).$ In particular, the stability conditions in Theorem \ref{mainthm} do not lie in the $\widetilde{\GL}_2^+(\R)$-orbit of Li's stability conditions.
\end{theorem}
    
\begin{proof}
    Suppose $Z^G=Z^G_{\omega,B,s,t}=M\circ Z^{\mathrm{Li}}_{a,b}$. Then,
    \[
    Z^G=x\Re Z^{\mathrm{Li}}+v\Im Z^{\mathrm{Li}}+i(u\Re Z^{\mathrm{Li}}+y\Im Z^{\mathrm{Li}}),\ \text{with $xy-uv>0$}.
    \]
    We now compare the coefficients of the independent functionals.
    \begin{enumerate}
        \item $\beta_2$ occurs only in $\Re Z^{\mathrm{Li}}$ with coefficient $-1/|G|$. It does not occur in $\Im Z^G.$ This implies $u=0.$ In $\Re Z^G$ it has coefficient $-1.$ So $x=|G|.$ Thus $y>0.$

        \item $\alpha_2$ and $\omega_0\cdot\beta_1$ occur in $Z^{\mathrm{Li}}$. The $\alpha_2$ coefficient of $\Re Z^G$ is $0.$ So $|G|(b+\tfrac{n+1}{2})+va=0.$ Then the $\beta_1$-functional on the right is $0.$ But on the left it is $\beta_1\mapsto B\cdot\beta_1,$ so $B=0.$ In the imaginary part: $\omega\cdot\beta_1=\tfrac{ya}{|G|}\omega_0\cdot\beta_1$ for all $\beta_1.$ This gives $\omega=\lambda\omega_0$ with $\lambda:=ya/|G|>0$. The $\alpha_2$-coefficients give $t=\lambda m.$ 

        \item $\alpha_1$ in the imaginary parts: the coefficient of the left side is $-tB=0$ by (2). On the right it is $-ya(b+\tfrac{n+1}{2})\cdot\tfrac{2m}{|G|}\omega_0.$ Hence $b=-\tfrac{n+1}{2},$ and then $v=0$ by (2).

        \item $\alpha_0$ in the imaginary part: With $B=0,\omega=\lambda\omega_0,$ and $b=-\tfrac{n+1}{2},$ we get,
        \begin{align*}
            \lambda m(1-\tfrac{\lambda^2\omega_0^2}{2})&=|G|\lambda K_2=\lambda[3m\omega_0^2(K_1-\tfrac{(n+1)^2}{8}-\tfrac{a^2}{6})+2m]\\
            \implies a^2&=6K_1-3(\tfrac{n+1}{2})^2+\lambda^2+\tfrac{2}{\omega_0^2}.
        \end{align*}

        \item $\beta_0$ in the real part: the coefficient on the left side is $-(1-\tfrac{\lambda^2\omega_0^2}{2}).$ The right side is $(-K_1+\tfrac{(n+1)^2}{8}+\tfrac{a^2}{2})\omega_0^2-2.$ So $a^2=2K_1-\tfrac{(n+1)^2}{4}+\lambda^2+\tfrac{2}{\omega_0^2}.$
    \end{enumerate}
    Notice, (4) and (5) are compatible only if $2K_1=(\tfrac{n+1}{2})^2$, i.e. $3n+4=3n+3$, hence a contradiction.
\end{proof}

\bibliographystyle{alpha}
\bibliography{mainbib}

@misc{LMPSZ,
 author = {Chunyi Li and Emanuele Macr{\`{\i}} and Alexander Perry and Paolo Stellari and Xiaolei Zhao},
 title = {Stability conditions on products of curves and {Hilbert} schemes of surfaces},
 year = {2025},
 howpublished = {Preprint, {arXiv}:2512.14207 [math.{AG}] (2025)},
 url = {https://arxiv.org/abs/2512.14207},
 arXiv = {arXiv:2512.14207}
}

@article{Liu-product,
 author = {Liu, Yucheng},
 title = {Stability conditions on product varieties},
 fjournal = {Journal f{\"u}r die Reine und Angewandte Mathematik},
 journal = {J. Reine Angew. Math.},
 issn = {0075-4102},
 volume = {770},
 pages = {135--157},
 year = {2021},
 language = {English},
 doi = {10.1515/crelle-2020-0010},
 zbMATH = {7291141},
 Zbl = {1468.14040}
}

@article{MMS,
 author = {Macr{\`{\i}}, Emanuele and Mehrotra, Sukhendu and Stellari, Paolo},
 title = {Inducing stability conditions},
 fjournal = {Journal of Algebraic Geometry},
 journal = {J. Algebr. Geom.},
 issn = {1056-3911},
 volume = {18},
 number = {4},
 pages = {605--649},
 year = {2009},
 language = {English},
 doi = {10.1090/S1056-3911-09-00524-4},
 zbMATH = {5610515},
 Zbl = {1175.14010}
}

@article{Bri-K3,
 author = {Bridgeland, Tom},
 title = {Stability conditions on {{\(K3\)}} surfaces},
 fjournal = {Duke Mathematical Journal},
 journal = {Duke Math. J.},
 issn = {0012-7094},
 volume = {141},
 number = {2},
 pages = {241--291},
 year = {2008},
 language = {English},
 doi = {10.1215/S0012-7094-08-14122-5},
 zbMATH = {5237682},
 Zbl = {1138.14022}
}

@misc{LLLMPSZ:26,
 author = {Chunyi Li and Zhiyu Liu and Ziqi Liu and Emanuele Macr{\`{\i}} and Alexander Perry and Paolo Stellari and Xiaolei Zhao},
 title = {Stability conditions and moduli spaces on projective families},
 year = {2026},
 howpublished = {Preprint, {arXiv}:2607.28411 [math.{AG}] (2026)},
 url = {https://arxiv.org/abs/2607.28411},
 arXiv = {arXiv:2607.28411}
}

@misc{Cheng:26,
 author = {Yiran Cheng},
 title = {A remark on the full support property},
 year = {2026},
 howpublished = {Preprint, {arXiv}:2608.14540 [math.{AG}] (2026)},
 url = {https://arxiv.org/abs/2608.14540},
 arXiv = {arXiv:2608.14540}
}

@article{BMS16,
 author = {Bayer, Arend and Macr{\`{\i}}, Emanuele and Stellari, Paolo},
 title = {The space of stability conditions on abelian threefolds, and on some {Calabi}-{Yau} threefolds},
 fjournal = {Inventiones Mathematicae},
 journal = {Invent. Math.},
 issn = {0020-9910},
 volume = {206},
 number = {3},
 pages = {869--933},
 year = {2016},
 language = {English},
 doi = {10.1007/s00222-016-0665-5},
 zbMATH = {6664764},
 Zbl = {1360.14057}
}

@article{BMT-I,
 author = {Bayer, Arend and Macr{\`{\i}}, Emanuele and Toda, Yukinobu},
 title = {Bridgeland stability conditions on threefolds. {I} {Bogomolov}-{Gieseker} type inequalities},
 fjournal = {Journal of Algebraic Geometry},
 journal = {J. Algebr. Geom.},
 issn = {1056-3911},
 volume = {23},
 number = {1},
 pages = {117--163},
 year = {2014},
 language = {English},
 doi = {10.1090/S1056-3911-2013-00617-7},
 zbMATH = {6276472},
 Zbl = {1306.14005}
}

@misc{Cheng-Feyzbaksh,
 author = {Yiran Cheng and Soheyla Feyzbakhsh},
 title = {Stability conditions on threefolds},
 year = {2026},
 howpublished = {Preprint, {arXiv}:2607.04788 [math.{AG}] (2026)},
 url = {https://arxiv.org/abs/2607.04788},
 arXiv = {arXiv:2607.04788}
}

@misc{FKLR-Brill-Noether-threefolds,
 author = {Soheyla Feyzbakhsh and Naoki Koseki and Zhiyu Liu and Nick Rekuski},
 title = {Stability conditions on {Calabi}-{Yau} threefolds via {Brill}-{Noether} theory of curves},
 year = {2025},
 howpublished = {Preprint, {arXiv}:2509.24990 [math.{AG}] (2025)},
 url = {https://arxiv.org/abs/2509.24990},
 arXiv = {arXiv:2509.24990}
}

@article{TypeK-classification,
 author = {Hashimoto, Kenji and Kanazawa, Atsushi},
 title = {Calabi-{Yau} three-folds of type {K}. {I}: {Classification}},
 fjournal = {IMRN. International Mathematics Research Notices},
 journal = {Int. Math. Res. Not.},
 issn = {1073-7928},
 volume = {2017},
 number = {21},
 pages = {6654--6693},
 year = {2017},
 language = {English},
 doi = {10.1093/imrn/rnw197},
 zbMATH = {7004483},
 Zbl = {1405.14098}
}

@misc{Perr-Shah,
 author = {Alexander Perry and Saket Shah},
 title = {Stability conditions on crepant resolutions of quotients of product varieties},
 year = {2024},
 howpublished = {Preprint, {arXiv}:2404.09121 [math.{AG}] (2024)},
 url = {https://arxiv.org/abs/2404.09121},
 arXiv = {arXiv:2404.09121}
}

@article{Families-paper,
 author = {Bayer, Arend and Lahoz, Mart{\'{\i}} and Macr{\`{\i}}, Emanuele and Nuer, Howard and Perry, Alexander and Stellari, Paolo},
 title = {Stability conditions in families},
 fjournal = {Publications Math{\'e}matiques},
 journal = {Publ. Math., Inst. Hautes {\'E}tud. Sci.},
 issn = {0073-8301},
 volume = {133},
 pages = {157--325},
 year = {2021},
 language = {English},
 doi = {10.1007/s10240-021-00124-6},
 zbMATH = {7395794},
 Zbl = {1481.14033}
}

@misc{Li,
 author = {Chunyi Li},
 title = {A {Remark} on {Stability} {Conditions} on {Smooth} {Projective} {Varieties}},
 year = {2026},
 howpublished = {Preprint, {arXiv}:2601.22994 [math.{AG}] (2026)},
 url = {https://arxiv.org/abs/2601.22994},
 arXiv = {arXiv:2601.22994}
}

@misc{LLPZ-Frobenius,
 author = {Chunyi Li and Zhiyu Liu and Alexander Perry and Xiaolei Zhao},
 title = {Generalized {Bogomolov}--{Gieseker} inequalities via {Frobenius}},
 year = {2026},
 howpublished = {Preprint, {arXiv}:2610.18595 [math.{AG}] (2026)},
 url = {https://arxiv.org/abs/2610.18595},
 arXiv = {arXiv:2610.18595}
}

@article{Oguiso-sakurai-type-A,
 author = {Oguiso, Keiji and Sakurai, Jun},
 title = {Calabi-{Yau} threefolds of quotient type},
 fjournal = {The Asian Journal of Mathematics},
 journal = {Asian J. Math.},
 issn = {1093-6106},
 volume = {5},
 number = {1},
 pages = {43--77},
 year = {2001},
 language = {English},
 doi = {10.4310/AJM.2001.v5.n1.a5},
 zbMATH = {1818524},
 Zbl = {1031.14022}
}

\end{document}